# History of confluent Vandermonde matrices and inverting them algorithms[1]

Jerzy S. Respondek[2]

*Silesian University of Technology*
*Applied Informatics Department*
*ul. Akademicka 16, 44-100 Gliwice, Poland*

**Abstract**

The author was encouraged to write this review by numerous requests from researchers around the world, who needed a ready-to-use algorithm for the inversion of confluent Vandermonde matrices which works in quadratic time for any values of the parameters allowed by the definition, including the case of large root multiplicity(-ties) of the characteristic polynomial. Article gives the history of the title matrix since 1891 and surveys algorithms for solving linear systems with the title class matrix and inverting it. It also presents a numerical algorithm which does not use symbolic computations and is ready to be implemented in a general-purpose programming language or in a mathematical package.

## 1. Introduction

Just as the classical Vandermonde matrix is associated with classical interpolation, its confluent generalisation is usually associated with generalised Hermite interpolation. However, in Section 4 we identified over a dozen other applications. Perhaps the most important of these is that the similarity matrix connecting the Jordan and Companion canonical forms of matrices is just the special matrix of the title.

In the literature there are many algorithms for its inversion and for solving systems of linear equations with this type of special matrix, but not all of them work for all input parameter series allowed by its definition. Moreover, some of them require symbolic computations, while others do not. Finally, only some of them are efficient in the general case.

In this survey, we have indicated those algorithms and methods that are free from the above mentioned drawbacks. The paper is structured as follows: in Section 2, we give a broad history of generic matrices. Using the examples of so-called compound matrices and determinoids, we also highlight the contemporary value of these seemingly historical and obsolete books. In Section 3 we have presented some other types of generalisation of the Vandermonde matrices, in addition to the title version. In Section 4 we have traced the history of the title matrices, enumerated over a dozen of their applications in different disciplines, and identified the main groups of algorithms for their inversion and for solving systems of equations with this type of special matrix. We have also proposed one which always works in quadratic time, and to work it needs only four basic arithmetic operations, without the need to perform any symbolic computations. In Section 5 we give some conclusions and perspectives. In the appendix, we give an example of the execution of this algorithm.

[1] This work is supported by the National Research Fund No 02/100/BK_25/0044

[2] jrespondek@polsl.pl, ORCID 0000-0001-5621-3783

## 2. An outline of the history of generic matrices

The concept of a two-dimensional mathematical matrix is a fundamental one in science and engineering. As documented, for example, in the monograph by Yong, Se [116], the concept of Gauss elimination, and thus an early version of matrices, was already known in ancient times. Section 7.7, pages 145–154, describes the process of transforming the source system matrix into a triangular form in order to solve a system of three linear equations with three unknowns.

Perhaps surprisingly, a much older concept than matrices are determinants. According to the 1890 monograph by Muir [82], they were invented as early as 1693 by Leibnitz, albeit in a different notation. For a comprehensive study of the history of determinants, see the classic and still valuable monographs by Muir [82]-[84]. Two notable monographs on determinants are Muir and Metzler [85] from 1933 and Scott [99] from 1880[3]. Another old book on determinants is Muir [81] from 1882. The concept of a matrix was first introduced by Sylvester in article [106] from 1850, on page 369, and again by Sylvester in article [107] from 1851, on page 302. Roughly speaking, Sylvester describes it as a creation, in general not necessarily square, from which we can pick a certain number of determinants. In other words, Sylvester defined the matrix as a subordinate concept to the determinants.

Nevertheless, many authors consider Cayley's 1857 article [16] to be the first comprehensive work devoted specifically to the topic of matrices. The 21-page text discusses 58 matrix operations and their properties. In addition to the widely recognised and appreciated monographs by Horn, Lancaster, Bellman and Gantmacher, it is also worth considering earlier works which are still valuable and treat matrix problems in a mature way. The above monographs are (we have deliberately presented them in reverse chronological order), Aitken [5] of 1944[4], preceded by Frazer, Duncan et al. [31] of 1938, Wedderburn [113] of 1934, Mac Duffee [26] of 1933[5], Turnbull & Aitken [110] of 1932[6], Turnbull [109] of 1928[7] and Cullis [23] vol. I-III dated in the range 1913-1925.

[3] A revised edition of this monograph by G.B. Mathews was published by the same publishers in 1904.

[4] This monograph went through several editions until the 9th edition in 1956, which was revised and expanded.

[5] A corrected reprint of the first edition of this monograph was published in 1956 by Chelsea Publishing Company, New York.

[6] The third edition of this book was published in 1952 and the last edition in 1961 by Dover Publications.

[7] The third edition of the book was published in 1960.

Interestingly, the monographs mentioned above, even the oldest ones, are not only of historical interest. An example is the concept of compound matrices. Let us first give their definition:

*2.1. Definition – compound matrices (e.g. Horn, Johnson (2013) [50] s.0.8.1 p. 21)*

*Let $A$ be a $m \times n$ dimensional matrix and $\alpha \subseteq \{1,...,m\}$ and $\beta \subseteq \{1,...,n\}$ be index sets of cardinality $r \le \min\{m,n\}$ elements. The $\binom{m}{r}$-by-$\binom{m}{r}$ matrix whose $\alpha,\beta$ entry is $\det A[\alpha,\beta]$ is called the $r$ th compound matrix of $A$ and is denoted by $C_r(A)$.*

An example of a compound matrix is given in the same section of the monograph [50]. Basic properties of compound matrices are presented in [50] p. 22. Compound matrices also retain a number of properties of the matrix on which they are built (also [50] p. 22).

These matrices are still used in many mathematical and engineering fields, including (in order of publication date):

- Network theory of electrical systems: Bryant (1963) [14].
- Tree-generating problem: Malik (1970) [77].
- Control theory:
  - General analysis of linear systems: Nambiar (1970) [87].
  - Controllability of control systems: The article Kalogeropoulos, Psarrakos (2004) [63] gave a new elegant and compact criterion for verifying the controllability of high order systems. It is worth noting that control theory, including controllability verification, was mature at the time of publication.
  - Bar-Shalom et al. (2023) [6] present an extensive tutorial on compound matrices in control theory.
- Matrix inequalities: articles Fiedler (1974) [29], Elsner et al. (2000) [27] used the compound matrices to prove a series of matrix inequalities involving their eigenvalues.
- Analysis of solid structures: Fuchs (1992) [32].
- Various differential equation problems:
  - Linsay, Rooney (1992) [67] (spectral problems of ordinary and partial differential equations).
  - Ivansson (2006) [56] (multi-point boundary-value problems).
  - Muldowney (1990) [86] (connection between compound matrices and ordinary differential equations).
- Graph theory: Article Nambiar (1997) [88] used compound matrices to give a concise and elegant proof of Hall's Theorem.
- Article Prells, Friswell (2003) [91] successfully tackled the classic hard problem of matrix algebra, i.e. finding the value of a determinant of the sum of matrices using compound matrices as a tool.
- Schur complements: Liu et al. (2010) [68].

Nevertheless, they are barely discussed (if at all) in contemporary monographs and textbooks. Gantmacher [35] devotes less than 4 pages to compound matrices in Section I.4 of Volume I. Even worse in this respect is the modern monograph, Horn, Johnson [50] from 2013. They describe compound matrices on a total of two pages, in sections 0.8.1 and 0.8.12.

On the other hand, a monograph by Cullis [23] of 1918 devotes an entire thirty-six page chapter XII to compound matrices. In Wedderburn [113] of 1934 we have a separate chapter V on pages 63-87. Probably the most cited classic book on matrices today, Aitken [5], devotes sections 38-41 and section 43 to compound matrices.

Another example of this kind is the generalization of determinants to rectangular matrices. Surprisingly, this generalization was already proposed in 1913 by Cullis [23], Vol. I, and called a '*determinoid*'[8]. It is defined there analogously to classical determinants, on a purely combinatorial basis. This concept did not reappear in the literature until 1966, when it was redefined on a higher but formally equivalent level by Radic [92]. In his definition he used the classical concept of determinants, which are appropriately linearly combined. The references of the article [76] include more than a dozen articles on the properties of determinoids and their specific applications, the first of which was published in 1991[9].

Determinoids cannot compete with compound matrices in the number of applications, but apparently this topic is being revived.

It is difficult to decide which of the contemporary monographs is the most comprehensive and reputable. But we decided to try. Probably the most widely known are: Horn, Johnson (2013) [50], Lancaster, Tismenetsky (1985) [64], Bellman (1970) [7] and Gantmacher (1960) [35][10]. A book specialized in matrix differential calculus is Magnus, Neudecker (2019) [75].

For a comprehensive contemporary monograph on the determinants, see Vein, Dale (1999) [111].

It is also worth mentioning that the monograph Golub, Van Loan [38] from 2013 is recognised as the most comprehensive work on numerical recipes for matrices. Another seminal work on numerical algorithms for matrices is Householder [54] from 1964. Another classic but still valuable book on numerical algorithms, focusing on the spectral properties of matrices, is Faddeev, Faddeeva (1963) [28]. Its earlier and shorter 1959 edition is also valuable (but only by the second author). A classic position dealing with various computational problems of matrices, but strongly related to error analysis, is Wilkinson (1988) [114].

To complete the argument, let us now give the definition of the matrix. In a monograph Bellman [7], page 16, the matrix is defined as a rectangular or, in particular, square table of complex numbers. Gantmacher [35], page 1, states that the entries of the 2D table must be elements of a certain field. Horn [50] and Lancaster [64] define it analogously. Some authors define the matrix concept in a more formal way. Kaczorek [57] on page 15 defines the matrix as "*the function of two variables which assigns to each pair of natural numbers (i,j), where i=1,..,m, j=1,..,n, exactly one real or complex number* $a_{ij}$". This means that Kaczorek defines it as a real or complex function defined on the Cartesian product $(1,..,m)\times(1,...,n)$.

[8] It is worth noting that the term *'determinoid'* appears in the very title of Cullis's monograph (1913).

[9] Not all generalisations of determinants to rectangular matrices are formally equivalent to the Cullis's definition.

[10] We give an English-translated edition of the Russian-language book "*Teorija Matrits*" by the same author.

Now we come to the operation of matrix multiplication. In point 11 of the article by Cayley [16], on pages 20-21, the definition of matrix multiplication in the case of a pair of matrices of dimension $3\times3$ is given (but in a notation very different from the contemporary one). In the early years of its development, matrix multiplication was also defined indirectly as an auxiliary tool for determinants. In his 1880 monograph [99], Scott gives an unintended definition of matrix multiplication in Chapter IV, entitled "*Multiplication of Determinants*". In section 3, pages 46–47, the defining formula is presented in the context of the formula for the determinant of a product of matrices. In monographs e.g. Gantmacher ([35], page 6), Bellman ([7], page 39), the operation is already defined in a contemporary way.

The matrix multiplication operation plays an important role because we have an extensive and quite mature separate branch of so-called fast matrix multiplication algorithms, of which only a tiny example is the milestone article by Coppersmith, Winograd [21]. There is also a mature branch of methods that allow to exploit the efficiency of fast matrix multiplication in other matrix and non-matrix operations, often with the same efficiency. Algorithms in this group owe their efficiency in part to exploiting the non-commutativity property of matrix multiplication. According to Hawkins [45], the non-commutativity of the matrix product was first identified by Eisenstein in works from 1844 and 1852 (see references therein). This fact is also stated by Cayley [16] in point 11.

## 3. Generalized Vandermonde matrices as a kind of special matrices

As we said in Section 2, one of the possible ways to make matrix operations faster are the so-called fast matrix multiplication algorithms. Their advantage is their generality, i.e. their input parameters are matrices with arbitrary allowed values for each entry. The disadvantage is that they are still quite time-consuming. Today, we do not know of any universal algorithm for matrix multiplication that is more efficient than $O(n^{2.37})$. Moreover, since 1990 (year of publication of the Coppersmith, Winograd [21] paper), the efficiency exponent has improved only by $\sim 0.003$. One can ask whether it is possible to improve their time complexity further.

On the other hand, further significant efficiency improvements are possible thanks to specialised algorithms designed for selected classes of special matrices. Special matrices are characterised by a specific structure of their entries. The following matrices can be distinguished: Pascal, Frobenius, Toeplitz, Hankel, Sylvester, Bézout, Cauchy, block diagonal, block triangular, companion, Jordan, circulant, Hilbert, bool matrices, just Vandermonde matrices (particularly generalised) and others. From the point of view of algorithm construction, their important property is that it is often possible to achieve $O(n^2)$ time complexity by the algorithms for certain parameter values of a given special matrix type, or even in the general case.

Unlike matrices of any form, it is difficult to provide a single, comprehensive monograph on special matrices. Commonly used types of matrices are Toeplitz and Hankel, with equal diagonal and anti-diagonal elements, respectively. They have been the subject of two specialised monographs: Iohvidov (1982) [55] and Heinig, Rost (1984) [47]. Böttcher,

Silbermann (2012) [11] deals with a kind of Toeplitz matrices. Davis (2012) [24] is devoted exclusively to circulant matrices. The multi-author books Bini et al. (2010) [8] and Kailath, Sayed (1999) [58] comprehensively present different types of numerical recipes for various types of structured matrices. Another monograph on fast algorithms that operate on various types of special matrices is Pan (2001) [89]. The article Aceto, Trigiante (2001) [1] is an appreciated work on this topic.

In the literature there are several directions of generalization of Vandermonde matrices, not only the confluent one. One of them is proposed by the monographs Gantmacher (1959) [34], Example 1 in Section 8, and Gantmacher (1960) [35], Example 1 also in Section 8 of Volume II, devoted to totally non-negative matrices[11], which is this generalisation:

$$V_{Gantmacher} = \begin{bmatrix} a_1^{\alpha_1} & a_2^{\alpha_1} & \cdots & a_n^{\alpha_1} \\ a_1^{\alpha_2} & a_2^{\alpha_2} & \cdots & a_n^{\alpha_2} \\ \vdots & \vdots & \ddots & \vdots \\ a_1^{\alpha_n} & a_2^{\alpha_n} & \cdots & a_n^{\alpha_n} \end{bmatrix}$$

for real parameters satisfying the inequalities:

$$\begin{aligned} 0 < a_1 < a_2 < \ldots < a_n \\ \alpha_1 < \alpha_2 < \ldots < \alpha_n \end{aligned} \tag{1}$$

The classical Vandermonde matrix is a special case of, let's call it, its Gantmacher generalisation, for $\alpha_i = i$; $i = 0,1,\ldots,n-1$. In Section 16, Vol. II, Gantmacher (1960) [35] applied this kind of Vandermonde generalisation to the problem of moments on the positive axis.

Probably the first paper in which this branch of the Vandermonde generalization was proposed is the article by Heineman (1929) [46]. The only small difference is that the exponents are assumed to be positive integers. He applied this matrix to certain algebraic constructions using elementary symmetric functions. Schlickewei, Viola (2000) [98] used this version of Vandermonde's generalization to solve recursive equations. Article Demmel, Koev (2005) [25] proposes an algorithm for solving linear equation set with this type of structured matrix, which is fast and accurate. The authors claim that the presented algorithm works on the virtue corresponding to that presented in the classical article Bjorck, Pereyra (1970) [10] adapted to the particular type of Vandermonde generalization (more details on the latter article we give in the next section). The article Chen et. al (2008) [20] gives a method for LU decomposition and recursive formulas for inversion and determinant value of this kind of generalized Vandermonde matrix.

[11] According to Section 8 in Gantmacher (1960) [35], Vol. II, a given matrix is totally nonnegative if and only if all its minors of any order are nonnegative.

Yet another kind of generalisation is proposed by El-Mikkawy, e.g. in article [78]:

$$V_{Mikkawy} = \begin{bmatrix} a_1^{k} & a_2^{k} & \cdots & a_n^{k} \\ a_1^{k+1} & a_2^{k+1} & \cdots & a_n^{k+1} \\ \vdots & \vdots & \ddots & \vdots \\ a_1^{k+(n-1)} & a_2^{k+(n-1)} & \cdots & a_n^{k+(n-1)} \end{bmatrix}$$

You can see that this is also a special case of the Gantmacher generalisation, now for $\alpha_i = k + i$ ; $i = 0,1,...,n-1$.

There is also a generalization that combines the classical (rectangular) Vandermonde matrix with another type of special matrix, the Cauchy matrix. It is constructed as a block matrix as follows (e.g., Finck, Heinig, Rost (1993) [30]):

$$\left[\begin{array}{cccc|cccc} \frac{1}{x_1-y_1} & \frac{1}{x_1-y_2} & \cdots & \frac{1}{x_1-y_k} & 1 & x_1 & \cdots & x_1^{n-k-1} \\ \frac{1}{x_2-y_1} & \frac{1}{x_2-y_2} & \cdots & \frac{1}{x_2-y_k} & 1 & x_2 & \cdots & x_2^{n-k-1} \\ \vdots & \vdots & \ddots & \vdots & \vdots & \vdots & \ddots & \vdots \\ \frac{1}{x_n-y_1} & \frac{1}{x_n-y_2} & \cdots & \frac{1}{x_n-y_k} & 1 & x_n & \cdots & x_n^{n-k-1} \end{array}\right] = \left[ C(\overline{x},\overline{y}) \mid V_{n-k}(\overline{x}) \right]$$

for pairwise distinct elements of $\overline{x}$ and $\overline{y}$ vectors. It is used in rational interpolation.

Calvetti, Reichel (1993) [15] propose to generalise a Vandermonde in a way that its entries are polynomials defined by a kind of three-term recurrence relation. The classical Vandermode matrix is a special case of this generalisation. Article [15] also gives an algorithm that inverts this generalization in quadratic time for any input values defined by this generalisation. We will refer to this article again in the next section when its further generalisation is described.

Finally, let's present the definition of the confluent Vandermonde matrix. In general, compared to the classical Vandermonde matrix, its confluent generalisation also contains derivatives of the columns multiplied by some constant factor.

*Due to Turnbull, Aitken (1932) [110], chapter VI page 60 or Kalman (1984) [60] pages 18-19, it is determined by the characteristic polynomial:*

$$p(s) = (s-\lambda_1)^{n_1}(s-\lambda_2)^{n_2}\cdot\ldots\cdot(s-\lambda_r)^{n_r} \tag{2}$$

*where* $\lambda_1, \lambda_2, ...\lambda_r$ *are given pair wise distinct real zeros with* $n_1 + ... + n_r = n$.

*The confluent Vandermonde matrix* $V$ *related to the zeros of* $p(s)$ *is defined to be the* $n \times n$ *block matrix* $V = [V_1 \mid V_2 \mid \cdots \mid V_r]$. *Let us denote as* $\overline{f}(\lambda)$ *the column* $\begin{bmatrix} 1 & \lambda & \cdots & \lambda^{n-1} \end{bmatrix}^T$ *and by* $\overline{f}^{(j)}(\lambda)$ *the* $j^{th}$ *derivative of this column. Now* $V_k$ *block is* $n \times n_k$ *matrix with columns* $\overline{f}^{(j)}(\lambda_k)/j!$ *for* $j=0,1,\ldots,n_j-1$. *We can also present the block in the form after differentiation as follows:*

$$\left(V_k\right)_{ij}=\begin{cases}\binom{i-1}{i-j}\lambda_k^{i-j}, & \text{for } i\geq j\\ 0, & \text{otherwise}\end{cases}$$

*for* $k=1,2,\ldots,r$ *;* $i=1,2,\ldots,n$ *and* $j=1,2,\ldots,n_k$.

## 4. Confluent Vandermonde matrices – history and progress of algorithms

Confluent Vandermonde matrices first appeared in 1901[12] in a form almost identical to that used today. The article Vogt [112] on page 350 gave an example of a determinant where one of the eigenvalues has a multiplicity equal to three. We can read there (in French): „*Ce déterminant est la **généralisation** de celui de **Vandermonde***". This determinant has been presented and reviewed in certain level on pages 183-184 of Muir 1930 [84].

In a classic monograph from 1923, Muir's „*History of Determinants*" [83], the confluent Vandermonde matrix is explicitly presented on pages 178-180 in exactly the form we use today, in a transposed determinant form, with decent dimensionality 6 by 6. Three different eigenvalues are considered, including a triple (by convention, the empty entries stand for zeros):

$$V=\left|\begin{array}{c:cc:ccc} 1 & 1 & & 1 & & \\ x & y & 1 & z & 1 & \\ x^2 & y^2 & 2y & z^2 & 2z & 1\\ x^3 & y^3 & 3y^2 & z^3 & 3z^2 & 3z\\ x^4 & y^4 & 4y^3 & z^4 & 4z^3 & 6z^2\\ x^5 & y^5 & 5y^4 & z^5 & 5z^4 & 10z^3 \end{array}\right|$$

The corresponding characteristic polynomial is:

$$p(s)=(s-x)(s-y)^2(s-z)^3$$

The next time the confluent Vandermonde matrix appears is in Turnbull's 1928 monograph [109], page 28, example 4. Here, however, we have a kind of regression compared to Muir in 1923 or even Vogt in 1901. It now appears as a determinant of rather timid parameters; we have an $4\times 4$ dimension with, what is even worse, a maximum multiplicity of entries equal to two.

The Vandermonde matrices, both in their classical and confluent versions, are mainly associated with their two main applications: in the so-called classical canonical form and in interpolation. Namely, the confluent Vandermonde matrix is the change-of-base matrix that

[12] Some authors consider the article Schendel (1891) [97] in German to be the first article in which this matrix appeared.

transforms the companion matrix into the Jordan canonical form, according to the matrix equation:

$$C \cdot V = V \cdot J \tag{3}$$

where $V$ is just a (generally confluent) Vandermonde matrix, $C$ is a companion matrix and $J$ is Jordan canonical form.

This fact was noted in two works in 1932. The first is the article Aitken [3], section 6 „*Relation between Two Forms of Canonical Matrix*”, pages 88-90. On page 89 we have the relation in question both in the form of a matrix equation and as an $4 \times 4$ example in formula (IV). The second work is a classic but still valuable monograph Turnbull, Aitken [110] Section VI.1:”*The Classical Canonical Form deduced from the Rational Form”*, pages 58-62. On page 60 we have an $6 \times 6$ example of this type of special matrix. For a more recent study of this idea, see also the 1984 article by Kalman [60].

Then the confluent Vandermonde matrix appears in the 1938 article by Aitken [4], formula 5 on page 288, as a tool for finding matrix eigenvalues[13].

In addition to the special matrix in question, a separate section VI.50 is devoted to it in the monograph Aitken [5] 1944, pages 119-121.

Furthermore, the determinant of the confluent Vandermonde matrix in the general case was found in 1932. A suitable formula can be found in the monograph Turnbull, Aitken [110], page 63, example 5 or in Aitken [5] from 1944, formula 7 on page 121.

As mentioned above, the classical Vandermonde matrix is often associated with polynomial interpolation. This is the case when we prescribe the values of the polynomial at given points. This type of interpolation was generalised by Hermite [48] in 1878, in an article written by this author in French. Namely, he considered not only the sole values of an interpolating polynomial, but also allowed to specify the values of the derivative of a polynomial of any order at any point in the series, starting their order from zero to a given one, with a step of one[14]. Nevertheless, Hermite formulated his generalised interpolation on the basis of the integral calculus in the complex domain. Aitken [5] gives an example of how this type of polynomial interpolation leads to the confluent Vandermonde matrix (formula 2, page 119, section 50). Next, the 1960 article by Spitzbart [104] formulated Hermite interpolation in a general case using this class of special matrices.

Another branch here is its generalisation to a 2D case, but already without the Vandermonde matrix[15]. Ahlin [2] from 1964 constructed a bivariate osculating interpolation polynomial $f(x, y)$, in which at each node both the values and its partial and mixed partial derivatives match prescribed values. However, in Ahlin [2] there is an additional restriction that at each node of a 2D grid we consider the same total order of derivatives or partial derivatives. This restriction removes the 1974 article by Chawla, Jayarajan [18].

[13] At the time of the article appeared, the eigenvalues were called ‘*latent roots*’.

[14] Without this additional assumption, we are faced with a Birkhoff interpolation, which is generally ambiguous. A monograph on Birkhoff interpolation is Lorentz, Jetter et al. [69].

[15] The Birkhoff interpolation in this case is the subject of a monograph [70].

In addition to the two classical applications discussed above, the confluent Vandermonde matrix is used in the following problems in mathematics, physics, engineering, and control theory:

- Computing the exponential of a matrix (Moler, V-Loan (2003) [79] - methods 11 and 13, Harris, Fillmore et al. (2001) [44]) and Luther, Rost (2004) [74].
- Computing matrix functions: Hermite interpolation for matrices (Higham (2008) [49], Section 1.2.2, especially Note 1.9 page 7).
- Calculation of the minimum polynomial of a matrix (Halidias (2024) [43]).
- Coding information in the Hermitian code (Lee K, O'Sullivan (2010) [65]).
- Optimization of the non-homogeneous differential equations (Gorecki (1968) [41]).
- Control theory:
  - Controllability analysis of linear systems (Ha, Gibson (1980) [42], Respondek (2008) [93]).
  - Quantised control of sampled-data systems (Shen, Tan et al. (2017) [100], Sun et al. (2020) [105]).
  - Stability analysis (Chen, Fu et al. (2017) [19]).
  - Properties of time-delay systems analysis (Boussaada, Niculescu (2016) [12][16]).
  - Root multiplicity analysis (Boussaada, Niculescu (2016) [13]).
- Time series analysis (Klein, Spreij (2003) [62]).
- Several classes of polynomials with repeating variables:
  - Schur polynomials (Serrano, Maximenko [39] (to appear, doi:10.1080/03081087.2025.2464639)).
  - A kind of symmetric polynomials (Serrano, Maximenko [40]).
- Discrete fractional Fourier transform (Moya-Cessa, Soto-Eguibar [80]).
- Linear recursive equations with constant coefficients: This class of equations can be transformed into a matrix form, with the companion matrix as the system matrix, and solved by finding a certain power of it. We find the power using the property (3) (Sousa et al. (2019) [103], Appendix A[17]).
- Vandermonde reduction of Bezoutians: G. Heinig and K. Rost in a chapter "*Introduction to Bezoutians*", pages 25-118 of the monograph Bini et al. (2010) [8] show that Bezoutians can be reduced to block diagonal form with the help of confluent Vandermonde matrices (Proposition 9.7 page 88).

[16] In this article also Birkhoff class of special matrices are broadly used.

[17] This method was first applied to the Fibonacci sequence in Silvester (1979) [101], which leads to $2\times 2$ matrices. But the case of multiple, in this case double, root is not considered. Article Kalman (1982) [59] generalised this method to equations of any order, but still limited to those with characteristic equation with single roots.

Let us now move on to the operations we want to perform on the Vandermonde matrix, i.e. inverting it and solving linear systems of equations with this matrix.

We have a variety of algorithms for **solving linear systems** with such a structure. Some of them work in quadratic (or even faster) time for any input parameter allowed by the definition of this matrix. Most of them are now classics written by famous authors. Let us analyse them, not necessarily in chronological order:

- As the first of them we can present the article Bjorck, Pereyra (1970) [10]. The basic idea used is the Newton interpolation, together with LU decomposition, which are applied progressively. The article gives a detailed example of how the proposed algorithm works for a single confluence node. However, the authors do not analyse the time complexity in the general case. This article inspired the two other algorithms, presented in articles [25] and [108], which are discussed in other paragraphs.
- An interesting generalisation of the previous Bjorck and Pereyra algorithm is proposed by Tang, Golub (1981) [108]. It presents an equivalent progressive algorithm to the corresponding one, generally using the same methods, but all computations are now formulated in block matrix form. This makes the article suitable for a parallel version of the algorithm. In particular, it can be used to solve linear systems with a classical Vandermonde matrix (without confluencies) in parallel. Unfortunately, the time complexity is only analysed for the non-confluent case.
- Another proposal for an algorithm that can efficiently deal with the solution of systems of linear equations with a confluent Vandermonde matrix is given by Galimberti, Pereyra (1971) [33]. Its working principle is to transform the source system into upper triangular form by a series of matrix pre-multiplications. The result is an algorithm that finds the vector of unknowns in quadratic time in the general case.
- Another article is proposed by Björck, Elfving (1973) [9]. Using the divided difference method, they derived an algorithm that works in quadratic time for confluent Vandermonde matrices in a general case.
- Definitely the **most modern papers** on the subject of solving systems of equations are the articles by Lu (1994, 1995) [71], [72]. They successfully combined divide and conquer techniques with fast Fourier transform, resulting in algorithms that **work in linear-polylogarithmic time**. Thus, they definitely present the fastest algorithm for the problem in question. It should be stressed that such high efficiency is only possible for the problem of solving a system of linear equations with a confluent Vandermonde system matrix. This will be explained in the next paragraph.

The problem of solving linear systems of equations (for any kind of system matrix) is a simpler problem than inverting the system matrix, in the sense that in the first problem we have $n$ unknowns, while in the second the number of unknowns is $n^2$. In the case of the matrix in question, this is manifested by the fact that after the Second World War a considerable number of algorithms for **inverting** the confluent Vandermonde matrix were published. However, unlike in the case of solving a system, none of them maintained efficiency for arbitrary values of the input parameters. We can distinguish the following main lines of research:

- **Algorithms that always work in quadratic time, but only in special cases.**

  The algorithm of Gautschi (1963) [36] only allows matrix parameters with multiplicities equal to two. If some of the eigenvalues have a higher multiplicity, the algorithm does not apply.

- **Algorithms that work for arbitrary values of input parameters, but are only efficient in special cases.**
  - The article by Zhong, Zhaoyong (1998) [117] further generalises the form of the Vandermonde matrix proposed in the article [15] mentioned in the previous section. Namely, he adds to it the confluencies, i.e. the differentiated columns multiplied by a certain constant factor. As the matrix proposed in article [15] reduces in a special case to the classical Vandermonde matrix, the matrix proposed by Zhong, Zhaoyong (1998) [117] reduces in a special case to the confluent Vandermonde matrix. The proposed algorithm for its inversion achieves quadratic efficiency only for small values of all eigenvalue multiplicities (compared to the matrix dimension).
  - The basic idea of the algorithm proposed in the classic article Luther, Rost (2004) [74] is to use the representation of the inversion by the transposition of the "source" confluent Vandermonde matrix to be inverted. Namely, in section 4 of [74] Luther, Rost show in Theorem 4.2 that the inversion is equal to the product of three matrices: The Hankel matrix constructed from the coefficient of the characteristic polynomial defining the matrix, the transposed confluent Vandermonde matrix to be inverted and a block Toeplitz matrix containing the coefficients of the partial fractional decomposition of the inverse of this characteristic polynomial (formula (4.2) therein on page 98).

    The confluent Vandermonde matrix is obtained in this article as a kind of by-product of solving the linear ordinary differential equation with constant coefficients by the matrix exponent method. The special matrix under consideration has a slightly different structure from the usual one. Namely, it is in a transposed form without the division of the columns by the term $j!$.

    Nevertheless, the constructed algorithms achieve quadratic efficiency only for small values of all eigenvalue multiplicities. Otherwise, the algorithms presented become cubic efficiency algorithms.
  - An earlier work with the same drawback is Kaufman [61] from 1969 (the author states on page 776, that it involves $k_i \times 5n^2$ operations, as $k_i$ denoted the i$^{th}$ root multiplicity). More precisely, Kaufman [61] also has embedded some symbolic computations, but they are easy to compute being of the form of finding derivatives of a characteristic polynomial (2). Yet another work of this kind is Chang [17] from 1974.

- **Inversion in the form of certain series of symbolic formulas to be differentiated.**

  There is a long line of such algorithms. Examples, in a chronological order, are: Schappelle [96] in 1972, Goknar [37] in 1973 and Csaki [22] in 1975.

- **Algorithms which returns only a part of the inversion entries.**

  The article Lupas [73] from 1975 is claimed to be fast, but the cost of that is for each block of the inversion (see further the block matrix formula for $W_i$ in the Theorem 4.1) gives only the last two rows (first paragraph on the page 560 of [73]).

- **Algorithms using the displacement structures**

  The notion of displacement structures is a powerful notion designed to operate on different classes of special matrices. For confluent Vandermonde matrices, the corresponding structure is given in the article by Yang et al. (2005) [115].

Until 2011, there was no algorithm in the literature for the inversion of the confluent Vandermonde matrix that simultaneously had the following capabilities:

- Works with all values of the input parameter series allowed by the definition of the particular special matrix.
- Always works in quadratic time, even for large multiplicities of one or a few eigenvalues.
- Maintains high accuracy for any input parameter series. This problem is presented in e.g. Pan [90] (2016) in the case of classical Vandemonde matrices and by Li [66] (2006) for the confluent generalisation of Vandemonde matrices.
- Does not require symbolic computation.

The breakthrough came only with the articles by Respondek [94], [95] (2011)[18], together with the previous work by Hou, Pang [53] (2002) (or Hou, Hou (2007) [51]), which met all four of the above requirements simultaneously. The steps in creating the algorithm through each paper are as follows:

- The first step in the construction of the inverse is made by the article Hou, Pang [53] (2002). It is based on the Leverrier-Faddeev construction of the resolvent of a Jordan matrix block (e.g. Faddeev, Faddeeva (1963) [28] Section 47 in Chapter IV, pp. 260-265) and then on a Taylor expansion. The confluent Vandermonde matrix itself appears after the second step, thanks to the differentiation in the nominator and the factorial in the denominator of the terms appearing in the Taylor series. The inversion is presented in a form that requires some symbolic calculations. A formally equivalent result is also given in Hou, Hou (2007) [51], where the confluent Vandermonde matrix is obtained by the similarity property (3) and the representation of the characteristic polynomial in a companion matrix form.

[18] The problem of solving a system of equations $Vx = b$ can be solved using matrix inversion by equality $x = V^{-1}b$.

- Further progress in the matter made by the article Respondek [94] is the simplification of the inversion formulas from the work of Hou, Pang [53] (2002) by analytical methods.
- The final step in the derivation of the algorithm is made in the article by Respondek [95]. He showed how a series of derivatives of $k=1,..,r$ different rational functions, derived from the characteristic polynomial (2) of the matrix, can be computed at once:
  - in a $k=1,..,r$ points simultaneously,
  - for any, but (in general) pairwise different, maximum order for each $k=1,..,r$, starting from 0 to $n_k-1$,

  using only the four basic arithmetical operations. The core of the algorithm is proved by an inductive method.

Overall, the works of Hou, Pang [53] (or Hou, Hou (2007) [51]) and Respondek [94], [95] lead to the following algorithm for inverting the confluent Vandermonde matrices:

*4.1. Algorithm – inversion of the confluent Vandermonde matrix*

The *inverse of the confluent Vandermonde matrix* $V$ *has the form* $V^{-1}=\left[W_1^T \mid W_2^T \mid \cdots \mid W_r^T\right]^T$. *The column vectors* $h_{k\cdot}$ *of the block matrix* $W_k=\left[h_{kn} \mid h_{k(n-1)} \mid \cdots \mid h_{k1}\right]$ *in the inverse confluent Vandermonde matrix* $V^{-1}$ *may be recursively computed by the following scheme:*

$$\begin{cases} h_{k1}=\left[K_{k,1}\cdots K_{k,n_k}\right]^T \\ h_{k2}=J_k\left(\lambda_k,n_k\right)h_{k1}+a_1h_{k1} \\ \vdots \\ h_{kn}=J_k\left(\lambda_k,n_k\right)h_{k(n-1)}+a_{n-1}h_{k1} \end{cases}, \quad k=1,2,...,r$$

where $a_k$ *are the coefficients of the characteristic polynomial* $p(s)=s^n+\sum_{i=0}^{n-1}a_{n-i}s^i$, $J_k\left(\lambda_k,n_k\right)$ *is the elementary Jordan block and* $K_{k,j}$ *are the auxiliary coefficients which may be computed by the following recursive scheme:*

$$\begin{cases} L_{ki}^{(q+1)}\left(\lambda_k\right)=q!\left(\lambda_k-\lambda_i\right)^q\cdot K_{k,n_k-q}-q\cdot L_{ki}^{(q)}\left(\lambda_k\right), \quad i=1,..,k-1,k+1,..,r \\ K_{k,n_k-q-1}=-\dfrac{1}{(q+1)!}\displaystyle\sum_{i=1,\,i\neq k}^{r}n_i\frac{L_{ki}^{(q+1)}\left(\lambda_k\right)}{\left(\lambda_k-\lambda_i\right)^{q+1}} \end{cases} \tag{4}$$

for $q=0,1,..,n_k-2$ *and* $k=1,2,...,r$. *The* $K_{k,n_k}$ *coefficients may be computed directly from the formula:*

$$K_{k,n_k}=\frac{1}{\left(\lambda_k-\lambda_1\right)^{n_1}\cdot\ldots\cdot\left(\lambda_k-\lambda_{k-1}\right)^{n_{k-1}}\left(\lambda_k-\lambda_{k+1}\right)^{n_{k+1}}\cdot\ldots\cdot\left(\lambda_k-\lambda_r\right)^{n_r}}, \quad k=1,2,...,r \tag{5}$$

The determination of the auxiliary coefficients $K_{k,j}$ by Theorem 4.1 deserves a closer look because of its sophisticated form. First we have to calculate the $K_{k,n_k}$ coefficient directly from the formula (5). Then, for each $k = 1, 2, ..., r$ and $q = 0, 1, .., n_k - 2$, we have to perform the following steps:

- Calculate the series $L_{k1}^{(q)}, L_{k2}^{(q)}, ..., L_{k(k-1)}^{(q)}, L_{k(k+1)}^{(q)}, ..., L_{kr}^{(q)}$ by the formula (4).
- Calculate the coefficient $K_{k,n_k-q-1}$ from the formula (4).

The order of the required calculations is detailed in Figure 1. The following should be noted:

- There is no need to initialise the $L_{ki}^{(0)}(\lambda_k)$, because for $q = 0$ they vanish, due to the formula (4).
- The calculation scheme shown in Figure 1 must be carried out separately for each $k = 1, 2, ..., r$. However, the order of the $k$ variable does not matter.
- In each iteration of the $q$ variable the previous $L_{ki}^{(q)}(\lambda_k)$ values may be overwritten by the new $L_{ki}^{(q+1)}(\lambda_k)$.

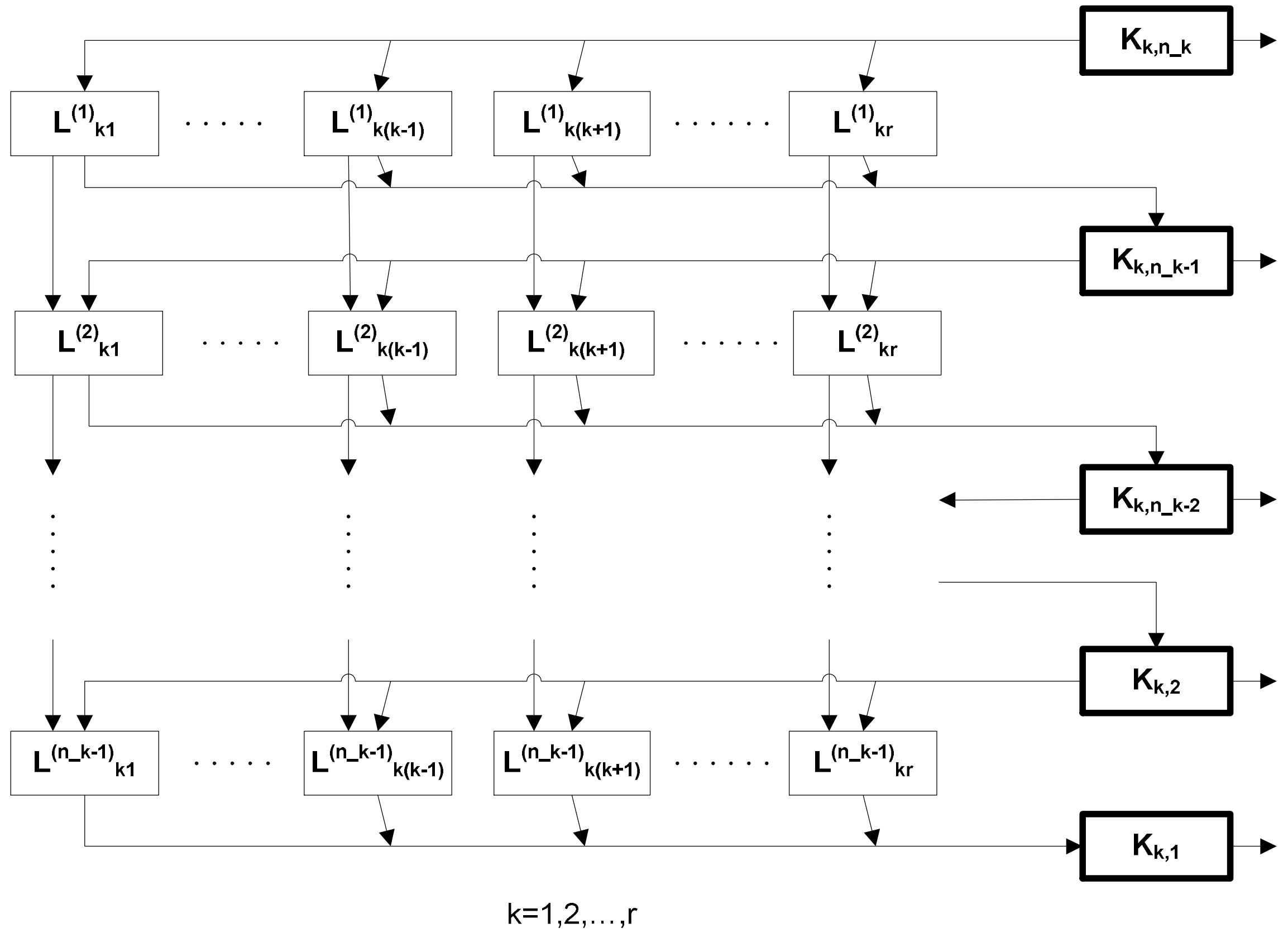

Figure 1 Scheme for the auxiliary coefficients determination algorithm

Theorem 4.1 allows to invert the confluent Vandermonde matrix without the need to perform symbolic computations. The article Respondek [95] provides a ready-to-use C/C++ implementation of the algorithm.

There are specialised algorithms for inverting a classical Vandermonde matrix that work in quadratic time, but for this purpose we can only use the algorithms for inverting its confluent version. Namely, the first one is a special case of the second for the $n_1 = n_2 = ... = n_r = 1$ parameters in formula (2).

To complete the lecture, let us show how to compute the determinant of the class of special matrices of the title. Turnbull, Aitken [110] (1932), page 63, example 5; or Aitken [5] (1944), formula 7 on page 121 for this purpose gives the following neat formula:

$$\det(V) = \prod_{1\le i<j\le r} \left(\lambda_j - \lambda_i\right)^{n_i n_j} \tag{6}$$

It is worth adding that the above formula is given in Aitken [5] and Turnbull, Aitken [110] without the detailed proof. The proof has been proposed in the article Sobczyk (2002) [102].

To summarise the possible representations of the inverse of the confluent Vandermonde matrix, we have available the following constructions:

- The triangular factorisation of the inverse, of the form $V^{-1} = HL$ - Hou, Hou [52] (2009).
- Matrix product combining the sought inverse with the transpose of the confluent Vandermonde matrix, of the form $V^{-1} = UV^T diag(DZP^T D)$, where $U$ is a Hankel matrix constructed from the coefficients of the characteristic polynomial defining the matrix, $V$ is the source generalized Vandermonde matrix to be inverted, $D$ is a block-diagonal matrix of constants, $Z$ is a block anti-diagonal matrix and $P$ is a block-Toeplitz matrix containing the coefficients of the partial fractional decomposition of the inverse of the characteristic polynomial - Luther, Rost (2004) [74].
- Matrix product containing transposed Vandermonde matrix and a Bezoutian - Bini et al. (2010) [8], Corollary 9.8 on page 88.
- Displacement structure - Yang et al. (2005) [115].

## 5. Conclusions and perspectives

In this article we have reviewed a number of algorithms for efficiently inverting confluent Vandermonde matrices and solving systems of linear equations with this class of system matrices. The described works present results both convenient for numerical purposes and in a form of symbolic formulas.

We think we have convinced the reader that sometimes it is worth looking at older literature, and we have made the reader aware of the different directions of generalisation of the Vandermonde matrix, its importance in applications, and what is the difference between solving systems of linear equations (SLE) and matrix inversion, and why SLE solving can be done faster than inversion. This is often confused in the literature and SLE algorithms are wrongly claimed to be the fastest, in a common group with algorithms inverting the title matrix.

Of course, the decision as to which algorithm to use should be made on a case-by-case basis for a particular application. The hardware capabilities should also be taken into account. For these reasons, it is not possible to judge which algorithm is the best.

In the view of the above, we hope that the decision on which algorithm to choose will be easier after reading this article.

Furthermore, the following can be considered as desirable future research directions in the field of numerical recipes for confluent Vandermonde matrices:

- Construction of the parallel algorithm for inverting the confluent Vandermonde matrices.
- Adaptation of the algorithm to the vector-oriented hardware units, like Intel AVX.
- Combination of both.

## References

[1] L. Aceto, D. Trigiante, The matrices of Pascal and other greats, Amer. Math. Monthly 108 (3) (2001) 232–245.

[2] A.C. Ahlin, A bivariate generalization of Hermite's interpolation formula, Mathematics of Computation 18 (1964) 264-273.

[3] A.C. Aitken, XII—Further Numerical Studies in Algebraic Equations and Matrices, Proc. Royal Society of Edinburgh 51 (1932) 80-90.

[4] A.C. Aitken, XX—Studies in Practical Mathematics. II. The Evaluation of the Latent Roots and Latent Vectors of a Matrix, Proc. Royal Society of Edinburgh 57 (1938) 269-304.

[5] A.C. Aitken, Determinants and Matrices, Third Edition, Interscience Publishers, New York, USA, 1944.

[6] E. Bar-Shalom, O. Dalin, M. Margaliot, Compound matrices in systems and control theory: a tutorial, Mathematics of Control, Signals, and Systems 35 (2023) 467–521.

[7] R. Bellman, Introduction to Matrix Analysis, 2nd. ed., McGraw-Hill, New York, USA, 1970.

[8] D.A. Bini, V. Mehrmann, V. Olshevsky, Numerical Methods for Structured Matrices and Applications, Operator Theory: Advances and Applications Book 199, Birkhäuser 2010.

[9] A. Björck, T. Elfving, Algorithms for confluent Vandermonde systems, Numerische Mathematik 21 (1973) 130-137.

[10] A. Bjorck, V. Pereyra, Solution of Vandermonde systems of equations, Math. Comp. 24 (1970) 893-903.

[11] A. Böttcher, B. Silbermann, Introduction to Large Truncated Toeplitz Matrices, Springer, 1st ed. 2012.

[12] I. Boussaada, SI Niculescu, Characterizing the Codimension of Zero Singularities for Time-Delay Systems. Acta Appl. Math. 145 (2016) 47–88.

[13] I. Boussaada, S.-I. Niculescu, Tracking the Algebraic Multiplicity of Crossing Imaginary Roots for Generic Quasipolynomials: A Vandermonde-Based Approach, IEEE Transactions on Automatic Control 61 (6) (2016) 1601-1606.

[14] P.R. Bryant, Compound Matrices In Network Theory, Recent Developments in Network Theory, Proc. Symposium Held at the College of Aeronautics, Cranfield, September 1961 (1963) 3-18.

[15] D. Calvetti, L. Reichel, Fast Inversion Of Vandermonde-Like Matrices Involving Orthogonal Polynomials, Bit Numerical Mathematics 33 (1993) 473-484.

[16] A. Cayley, A Memoir on the Theory of Matrices, Philosophical Transactions of the Royal Society of London 148 (1857) 17-37.

[17] F-Ch. Chang, The inverse of the generalized Vandermonde matrix through the partial fraction expansion, IEEE Transactions on Automatic Control 19 (2) (1974) 151-152.

[18] M.M. Chawla, N. Jayarajan, A Generalization of Hermite's Interpolation Formula in two Variables, Journal of the Australian Mathematical Society 18 (4) (1974) 402-410.

[19] J. Chen, P. Fu, C.-F. Méndez-Barrios, S.-I. Niculescu, H. Zhang, Stability Analysis of Polynomially Dependent Systems by Eigenvalue Perturbation, IEEE Transactions on Automatic Control, 62 (1) (2017) 5915-5922.

[20] Y.-M. Chen, H.-Ch. Li, E.-T. Tan, An Explicit Factorization Of Totally Positive Generalized Vandermonde Matrices Avoiding Schur Functions, Applied Mathematics E-Notes 8 (2008) 138-147.

[21] D. Coppersmith, S. Winograd, Matrix multiplication via Arithmetic Progressions, J. Symbolic Computation 9 (1990) 251-280.

[22] F. Csaki, Some notes on the inversion of confluent Vandermonde matrices, IEEE Transactions on Automatic Control 20 (1) (1975) 154-157.

[23] C.E. Cullis, Matrices and Determinoids, Vol. 1-3, Cambridge University Press, Cambridge, UK, 1913, 1918, 1925.

[24] P.J. Davis, Circulant Matrices, American Mathematical Society 338, Chelsea Pub Co, 2nd ed. 2012.

[25] J. Demmel, P. Koev, The Accurate and Efficient Solution of a Totally Positive Generalized Vandermonde Linear System, SIAM J. Matrix Anal. Appl. 27 (1) (2005) 142–152.

[26] C.C. Mac Duffee, The Theory of Matrices, Springer Verlag, Berlin, Germany, 1933.

[27] L. Elsner, D. Hershkowitz, H Schneider, Bounds On Norms of Compound Matrices and on Products of Eigenvalues, Bull. London Math. Soc. 32 (2000) 15-24.

[28] V.N. Faddeev, D.K. Faddeeva, Computational Methods of Linear Algebra, W. h. Freeman Company 1963.

[29] M. Fiedler, Additive compound matrices and an inequality for eigenvalues of symmetric stochastic matrices, Czechoslovak Mathematical Journal 24 (3) (1974) 392–402.

[30] T. Finck, G. Heinig, K. Rost, An Inversion Formula and Fast Algorithms for Cauchy-Vandermonde Matrices, Linear Algebra and Its Applications 183 (1993) 179-191.

[31] R.A. Frazer, W.J. Duncan, A.R. Collar, Elementary matrices and some applications to dynamics and differential equations, Cambridge University Press, Cambridge, UK, 1938.

[32] M.B. Fuchs, The explicit inverse of the stiffness matrix, Int. J. Solids Struct. 29 (1992) 2101-2113.

[33] G. Galimberti, V. Pereyra, Solving confluent Vandermonde systems of Hermite type, Numerische Mathematik 18 (1971) 44–60.

[34] F.R. Gantmacher, Applications of the Theory of Matrices, Interscience Publishers, Inc., New York, USA, 1959.

[35] F.R. Gantmacher, The Theory of Matrices, Vol. I-II, Chelsea Publishing Company, New York, USA, 1960.

[36] W. Gautschi, On inverses of Vandermonde and confluent Vandermonde matrices. II, Numerische Mathematik 5 (1963) 425-430.

[37] I. Goknar, Obtaining the inverse of the generalized Vandermonde matrix of the most general type, IEEE Transactions on Automatic Control 18 (5) (1973) 530-532.

[38] G.H. Golub, Ch.F. Van Loan, Matrix Computations, Fourth Edition, The Johns Hopkins University Press, Baltimore, USA, 2013.

[39] L.A. González-Serrano, E.A. Maximenko, Bialternant formula for Schur polynomials with repeating variables, Linear and Multilinear Algebra, 2025 (accepted at doi:10.1080/03081087.2025.2464639).

[40] L.A. González-Serrano, E.A. Maximenko, Complete Homogeneous Symmetric Polynomials with Repeating Variables, Mathematics 13 (1) (2025) 34.

[41] H. Gorecki, On switching instants in minimum-time control problem, One-dimensional case n-tuple eigenvalue, Bull. de L'Acad. Pol. Des. Sci. 16 (1968) 23–30.

[42] T.T. Ha, J.A. Gibson, A note on the determinant of a functional confluent Vandermonde matrix and controllability, Linear Algebra and its Applications 30 (1980) 69-75.

[43] N. Halidias, Computing the Minimum Polynomial, the Function and the Drazin Inverse of a Matrix with Matlab, Asian Journal of Research in Computer Science 17 (5) (2024) 1-9.

[44] W.A. Harris, J.P. Fillmore, D.R. Smith, Matrix Exponentials-Another Approach, SIAM Review 43 (4) (2001) 694-706.

[45] T. Hawkins, The Theory of Matrices in the 19th Century, Proc. Int. Congress of Mathematicians, Vancouver, Vol. 2 (1974) 561-570.

[46] E.R. Heineman, Generalized Vandermonde Determinants, Trans. Amer. Math. Soc. 31 (1929) 464-476.

[47] G. Heinig, K. Rost, Algebraic Methods for Toeplitz-like Matrices and Operators, Springer Basel AG 1984.

[48] M.Ch. Hermite, M. Borchardt, Sur la formule d'interpolation de Lagrange, Journal für die reine und angewandte Mathematik 84 (1878) 70-79.

[49] N.J. Higham, Functions of Matrices, Theory and Computation, SIAM, Philadelphia, USA, 2008.

[50] R.A. Horn, Ch.R. Johnson, Matrix Analysis, Second Edition, Cambridge University Press, New York, USA, 2013.

[51] S.-H. Hou, E. Hou, Recursive Computation of Inverses of Confluent Vandermonde Matrices, Electronic J. Mathematics and Technology, 1 (1) (2007) 11-24.

[52] S-H. Hou, E.S-H. Hou, A Recursive Algorithm for Triangular Factorization of Inverse of Confluent Vandermonde Matrices, AIP Conf. Proc. Vol. 1089 (1) (2009) 277–288.

[53] S-H. Hou, W-K. Pang, Inversion of confluent Vandermonde matrices, Computers & Mathematics with Applications 43 (12) (2002) 1539-1547.

[54] A.S. Householder, The Theory of Matrices in Numerical Analysis, First edition, Blaisdell Publishing Company, New York, USA, 1964.

[55] I.S. Iohvidov, Hankel and Toeplitz Matrices and Forms, Algebraic Theory, Birkhauser, Boston 1982.

[56] S. Ivansson, The Compound Matrix Method for Multi-Point Boundary-Value Problems, Journal of Applied Mathematics and Mechanics 77 (10) (2006) 767-776.

[57] T. Kaczorek, Vectors and matrices in automation and electrical engineering, 2nd ed., WNT, Warsaw, Poland, 1998 (in Polish).

[58] T. Kailath, A.H. Sayed, Fast Reliable Algorithms for Matrices with Structure, SIAM 1999.

[59] D. Kalman, Generalized Fibonacci Numbers by Matrix Methods. The Fibonacci Quarterly 20 (1) (1982) 73–76.

[60] D. Kalman, The Generalized Vandermonde Matrix, Mathematics Magazine 57 (1) (1984) 15-21.

[61] I. Kaufman, The inversion of the Vandermonde matrix and transformation to the Jordan canonical form, IEEE Transactions on Automatic Control 14 (6) (1969) 774-777.

[62] A. Klein, P. Spreij, Some Results on Vandermonde Matrices with an Application to Time Series Analysis, SIAM J. Matrix Anal. Appl. 25 (1) (2003) 213–223.

[63] G. Kalogeropoulos, P. Psarrakos, A note on the controllability of higher-order linear systems, Applied Mathematics Letters 17 (12) (2004) 1375-1380.

[64] P. Lancaster, M. Tismenetsky, The Theory of Matrices: Second Edition With Applications, Computer Science and Applied Mathematics, Academic Press, 2nd ed. 1985.

[65] K. Lee, M.E. O'Sullivan, Algebraic soft-decision decoding of Hermitian codes, IEEE Trans. Inform. Theory 56 (6) (2010) 2587–2600.

[66] R-C. Li, Lower Bounds for the Condition Number of a Real Confluent Vandermonde Matrix, Mathematics of Computation 75 (256) (2006) 1987–1995.

[67] K.A. Linsay, C.E. Rooney, A note on compound matrices, J. Computat. Phys. 103 (1992) 472-477.

[68] J. Liu, R. Huang, Generalized Schur Complements of Matrices and Compound Matrices, Electronic Journal of Linear Algebra 21 (2010) 12-24.

[69] G.G. Lorentz, K. Jetter, S.D. Riemenschneider, Birkhoff Interpolation, 1st ed., Encyclopedia of Mathematics and its Applications Book 19, Cambridge University Press, UK, 1984.

[70] R.A. Lorentz, Multivariate Birkhoff Interpolation, Lecture Notes in Mathematics, Springer-Verlag 1516 (1992).

[71] H. Lu, Fast Solution of Confluent Vandermonde Linear Systems, SIAM J. Matrix Anal. Appl. 15 (4) (1994) 1277-1289.

[72] H. Lu, Fast Algorithms For Confluent Vandermonde Linear Systems and Generalized Trummers Problem, SIAM J. Matrix Anal. Appl. 16 (2) (1995) 655-674.

[73] L. Lupas, On the computation of the generalized Vandermonde matrix inverse, IEEE Transactions on Automatic Control 20 (4) (1975) 559-561.

[74] U. Luther, K. Rost, Matrix exponentials and inversion of confluent Vandermonde matrices, Electronic Transactions on Numerical Analysis 18 (2004) 91-100.

[75] J.R. Magnus, H. Neudecker, Matrix Differential Calculus with Applications in Statistics and Econometrics, Wiley Series in Probability and Statistics, Wiley, 3rd edition 2019.

[76] A. Makarewicz, P. Pikuta, Cullis–Radic determinant of a rectangular matrix which has a number of identical columns, Annales Universitatis Mariae Curie-Sklodowska LXXIV (2) (2020) 41–60.

[77] R.N. Malik, Compound matrices to the tree-generating problem, IEEE Trans. Circuit Theory 17 (1970) 149-151.

[78] M.E.A. El-Mikkawy, Explicit inverse of a generalized Vandermonde matrix, Applied Mathematics and Computation 146 (2003) 643–651.

[79] C. Moler, Ch. V-Loan, Nineteen Dubious Ways to Compute the Exponential of a Matrix, Twenty-Five Years Later, SIAM Review 45 (1) (2003) 3-49.

[80] H.M. Moya-Cessa, F. Soto-Eguibar, Discrete fractional Fourier transform: Vandermonde approach, IMA Journal of Applied Mathematics 83 (6) (2018) 908–916.

[81] T. Muir, A treatise on the theory of determinants, MacMillan and co., London, UK, 1882.

[82] T. Muir, The Theory of Determinants in the Historical Order of Its Development. Part I. Determinants in General. Leibnitz (1693) to Cayley (1841), MacMillan and co., London, UK, 1890.

[83] T. Muir, The Theory of Determinants in the Historical Order of Development. Volumes Two to Four. The Periods 1841 to 1860, 1861 to 1880, 1880 to 1900, Dover Publications, Inc., New York, USA, 1911, 1920, 1923.

[84] T. Muir, Contributions to the History of Determinants. 1900-1920, Blackie & Son Limited, London and Glasgow, UK, 1930.

[85] T. Muir, W.H. Metzler, Theory of Determinants, Longmans, Green and co., Bombay, Calcutta, Madras, India, 1933.

[86] J.S. Muldowney, Compound Matrices and Ordinary Differential Equations, The Rocky Mountain Journal of Mathematics 20 (4) (1990) 857-872.

[87] K.K. Nambiar, J.D Keating, Application of compound matrices to linear systems, IEEE Trans. Circuit Theory 17 (1970) 626-628.

[88] K.K. Nambiar, Hall's theorem and compound matrices, Math. Comput. Modelling 25 (1997) 23-24.

[89] V. Pan, Structured Matrices and Polynomials Unified Superfast Algorithms, Springer Science+Business Media, LLC, New York, USA, 2001.

[90] V. Pan, How Bad Are Vandermonde Matrices? SIAM J. Matrix Anal. Appl. 37 (2) (2016) 676-694.

[91] U. Prells, M.I. Friswell, S.D. Garvey, Use of geometric algebra: compound matrices and the determinant of the sum of two matrices, Proc. Royal Society 459 (2003) 273-285.

[92] M. Radic, A definition of determinant of rectangular matrix, Glas. Mat. Ser. III 1 (21) (1966) 17–22.

[93] J.S. Respondek, Approximate controllability of the n-th order infinite dimensional systems with controls delayed by the control devices, Int. J. Syst. Sci. 39 (8) (2008) 765–782.

[94] J.S. Respondek, On the confluent Vandermonde matrix calculation algorithm, Appl. Math. Lett. 24 (2011) 103–106.

[95] J.S. Respondek, Numerical recipes for the high efficient inverse of the confluent Vandermonde matrices, Appl. Math. Comput. 218 (2011) 2044–2054.

[96] R. Schappelle, The inverse of the confluent Vandermonde matrix, IEEE Transactions on Automatic Control 17 (5) (1972) 724-725.

[97] L. Schendel, Das alternirende Exponentialdifferenzenproduct, Zeitschrift Math. Phys. (1891) 84-94.

[98] H.P. Schlickewei, C. Viola, Generalized Vandermonde determinants, Acta Arithmetica XCV 2 (2000) 123-137.

[99] R.F. Scott, Theory of Determinants, Cambridge University Press, Cambridge, UK, 1880.

[100] B. Shen, H. Tan, Z. Wang, T. Huang, Quantized/Saturated Control for Sampled-Data Systems Under Noisy Sampling Intervals: A Confluent Vandermonde Matrix Approach, IEEE Transactions on Automatic Control 62 (9) (2017) 4753-4759.

[101] J.R. Silvester, Fibonacci Properties by Matrix Methods, The Mathematical Gazette 63 (425) (1979) 188-191.

[102] G. Sobczyk, Generalized Vandermonde determinants and applications, Aportaciones Matematicas, Serie Comunicaciones 30 (2002) 41-53.

[103] Y.R. Sousa, H. Takayasu, M. Takayasu, Random coefficient autoregressive processes and the PUCK model with fluctuating potential, Journal of Statistical Mechanics 1 (2019) 013403.

[104] A. Spitzbart, A Generalization of Hermite's Interpolation Formula, The American Mathematical Monthly 67 (1) (1960) 42-46.

[105] H. Sun, J. Sun, J. Chen, LQG control for sampled-data systems under stochastic sampling, Journal of the Franklin Institute, 357 (5) (2020) 2773-2790.

[106] J.J. Sylvester, Additions to the articles in the September number of this journal, "On a new class of theorems," and on Pascal's theorem, Phil. Mag. 3 (37) (1850) 363-370.

[107] J.J. Sylvester, On the Relation between the Minor Determinants of Linearly Equivalent Quadratic Functions, Phil. Mag. 4 (1) (1851) 295-305.

[108] W.P. Tang, G.H. Golub, The block decomposition of a Vandermonde matrix and its applications, BIT Numerical Mathematics 21 (1981) 505-517.

[109] H.W. Turnbull, The Theory of Determinants, Matrices and Invariants, Blackie & Son, London & Glasgow, UK, 1928.

[110] H.W. Turnbull, A.C. Aitken, An Introduction to the Theory of Canonical Matrices, London, Glasgow and Bombay: Blackie and Son, 1932.

[111] R. Vein, P. Dale, Determinants and Their Applications in Mathematical Physics, Applied Mathematical Sciences 134, Springer 1999.

[112] M. Vogt, Sur l'apolarité des formes binaires. Nouvelles annales de mathématiques 4 (1) (1901) 337-365.

[113] J.H.M. Wedderburn, Lectures on Matrices. American Mathematical Society, Colloquium Publications, Providence, Rhode Island, USA, 1934.

[114] J.H. Wilkinson, The Algebraic Eigenvalue Problem, Numerical Mathematics and Scientific Computation, Clarendon Press, Revised ed. 1988.

[115] Z. Yang, L. Wang, Y. Hu, Displacement structures and fast inversion formulas for confluent polynomial Vandermonde-like matrices, J. Computational and Applied Mathematics 180 (2) (2005) 229-243.

[116] L.L. Yong, A.T. Se, Fleeting Footsteps. Tracing the Conception of Arithmetic and Algebra in Ancient China, World Scientific Publishing, Singapore, 2004.

[117] X. Zhong, Y. Zhaoyong, A Fast Algorithm for Inversion of Confluent Vandermonde-Like Matrices Involving Polynomials That Satisfy a Three-Term Recurrence Relation, SIAM J. Matrix Anal. Appl. 19 (3) (1998) 797-806.

## Appendix – example of execution of the algorithm

Let us invert the matrix defined by the following characteristic polynomial:

$$p(s) = (s+0.5)(s+3)^2(s+2)^3(s+1)^4$$

We can rewrite its data in a tabular form:

**Table 1 Parameters of confluent Vandermonde matrix to invert**

| i | 1 | 2 | 3 | 4 |
|---|---|---|---|---|
| $\lambda_i$ | -0.5 | -3.0 | -2.0 | -1.0 |
| $n_i$ | 1 | 2 | 3 | 4 |

The expansion of the $p(s)$ polynomial can be done in quadratic time at worst. We can receive:

**Table 2 Characteristic polynomial coefficients after expansion**

| $i$<br>$i$ | 1 | 2 | 3 | 4 | 5 | 6 | 7 | 8 | 9 | 10 |
|---|---|---|---|---|---|---|---|---|---|---|
| $a_i$<br>$a_i$ | 16.5 | 119.0 | 493.5 | 1302.0 | 2281.5 | 2687.0 | 2098.5 | 1039.0 | 294.0 | 36.0 |

The matrix to be inverted has the following form:

$$V = [V_1 \mid V_2 \mid V_3 \mid V_4] = \left[\begin{array}{c|cc|ccc|cccc}
1 & 1 & 0 & 1 & 0 & 0 & 1 & 0 & 0 & 0 \\
-0.5 & -3 & 1 & -2 & 1 & 0 & -1 & 1 & 0 & 0 \\
0.25 & 9 & -6 & 4 & -4 & 1 & 1 & -2 & 1 & 0 \\
-0.125 & -27 & 27 & -8 & 12 & -6 & -1 & 3 & -3 & 1 \\
0.0625 & 81 & -108 & 16 & -32 & 24 & 1 & -4 & 6 & -4 \\
-0.0312 & -243 & 405 & -32 & 8 & -80 & -1 & 5 & -10 & 10 \\
0.01562 & 729 & -1458 & 64 & -192 & 240 & 1 & -6 & 15 & -20 \\
-0.0078 & -2187 & 5103 & -128 & 448 & -672 & -1 & 7 & -21 & 35 \\
0.0039 & 6561 & -17496 & 256 & -1024 & 1792 & 1 & -8 & 28 & -56 \\
-0.0019 & -19683 & 59049 & -512 & 2304 & -4608 & -1 & 9 & -36 & 84
\end{array}\right]$$

Table 3 shows the execution of Theorem 4.1.

**Table 3 Consecutive iterations of Theorem 4.1**

| $k$ | $\lambda_k$ | $n_k$ | $q$ | $i$ | $L_{ki}^{(q+1)}(\lambda_k)$ | $K_{k,n_k-q-1}$ |
|---|---|---|---|---|---|---|
| 1 | -0.5 | 1 | - | - | - | $K_{1,1}=0.75851$ |
| 2 | -3 | 2 | - | - | - | $K_{2,2}=0.02500$ |
| | | | 0 | 1 | $L_{21}^{(1)}(\lambda_2)=0.02500$ | $K_{2,1}=0.13500$ |
| | | | | 3 | $L_{23}^{(1)}(\lambda_2)=0.02500$ | |
| | | | | 4 | $L_{24}^{(1)}(\lambda_2)=0.02500$ | |
| 3 | -2 | 3 | - | - | - | $K_{3,3}=-0.66666$ |
| | | | 0 | 1 | $L_{31}^{(1)}(\lambda_3)=-0.66666$ | $K_{3,2}=-1.77777$ |
| | | | | 2 | $L_{32}^{(1)}(\lambda_3)=-0.66666$ | |
| | | | | 4 | $L_{34}^{(1)}(\lambda_3)=-0.66666$ | |
| | | | 1 | 1 | $L_{31}^{(2)}(\lambda_3)=3.33333$ | $K_{3,1}=-4.51851$ |
| | | | | 2 | $L_{32}^{(2)}(\lambda_3)=-1.11111$ | |
| | | | | 4 | $L_{34}^{(2)}(\lambda_3)=2.44444$ | |
| 4 | -1 | 4 | - | - | - | $K_{4,4}=-0.50000$ |
| | | | 0 | 1 | $L_{41}^{(1)}(\lambda_4)=-0.50000$ | $K_{4,3}=1.00000$ |
| | | | | 2 | $L_{42}^{(1)}(\lambda_4)=-0.50000$ | |
| | | | | 3 | $L_{43}^{(1)}(\lambda_4)=-0.50000$ | |
| | | | 1 | 1 | $L_{41}^{(2)}(\lambda_4)=0.00000$ | $K_{4,2}=-2.87500$ |
| | | | | 2 | $L_{42}^{(2)}(\lambda_4)=2.50000$ | |
| | | | | 3 | $L_{43}^{(2)}(\lambda_4)=1.50000$ | |
| | | | 2 | 1 | $L_{41}^{(3)}(\lambda_4)=-1.43750$ | $K_{4,1}=3.62500$ |
| | | | | 2 | $L_{42}^{(3)}(\lambda_4)=-28.0000$ | |
| | | | | 3 | $L_{43}^{(3)}(\lambda_4)=-8.75000$ | |

We can now obtain the following inversion of the confluent Vandermonde matrix:

$$
V^{-1}=\left[\begin{array}{c} W_1 \\ \hline W_2 \\ \hline W_3 \\ \hline W_4 \end{array}\right]=\left[\begin{array}{ccccccccc|c}
h_{1(10)}^{(1)} & h_{19}^{(1)} & h_{18}^{(1)} & h_{17}^{(1)} & h_{16}^{(1)} & h_{15}^{(1)} & h_{14}^{(1)} & h_{13}^{(1)} & h_{12}^{(1)} & K_{1,1} \\ \hline
h_{2(10)}^{(1)} & h_{29}^{(1)} & h_{28}^{(1)} & h_{27}^{(1)} & h_{26}^{(1)} & h_{25}^{(1)} & h_{24}^{(1)} & h_{23}^{(1)} & h_{22}^{(1)} & K_{2,1} \\
h_{2(10)}^{(2)} & h_{29}^{(2)} & h_{28}^{(2)} & h_{27}^{(2)} & h_{26}^{(2)} & h_{25}^{(2)} & h_{24}^{(2)} & h_{23}^{(2)} & h_{22}^{(2)} & K_{2,2} \\ \hline
h_{3(10)}^{(1)} & h_{39}^{(1)} & h_{38}^{(1)} & h_{37}^{(1)} & h_{36}^{(1)} & h_{35}^{(1)} & h_{34}^{(1)} & h_{33}^{(1)} & h_{32}^{(1)} & K_{3,1} \\
h_{3(10)}^{(2)} & h_{39}^{(2)} & h_{38}^{(2)} & h_{37}^{(2)} & h_{36}^{(2)} & h_{35}^{(2)} & h_{34}^{(2)} & h_{33}^{(2)} & h_{32}^{(2)} & K_{3,2} \\
h_{3(10)}^{(3)} & h_{39}^{(3)} & h_{38}^{(3)} & h_{37}^{(3)} & h_{36}^{(3)} & h_{35}^{(3)} & h_{34}^{(3)} & h_{33}^{(3)} & h_{32}^{(3)} & K_{3,3} \\ \hline
h_{4(10)}^{(1)} & h_{49}^{(1)} & h_{48}^{(1)} & h_{47}^{(1)} & h_{46}^{(1)} & h_{45}^{(1)} & h_{44}^{(1)} & h_{43}^{(1)} & h_{42}^{(1)} & K_{4,1} \\
h_{4(10)}^{(2)} & h_{49}^{(2)} & h_{48}^{(2)} & h_{47}^{(2)} & h_{46}^{(2)} & h_{45}^{(2)} & h_{44}^{(2)} & h_{43}^{(2)} & h_{42}^{(2)} & K_{4,2} \\
h_{4(10)}^{(3)} & h_{49}^{(3)} & h_{48}^{(3)} & h_{47}^{(3)} & h_{46}^{(3)} & h_{45}^{(3)} & h_{44}^{(3)} & h_{43}^{(3)} & h_{42}^{(3)} & K_{4,3} \\
h_{4(10)}^{(4)} & h_{49}^{(4)} & h_{48}^{(4)} & h_{47}^{(4)} & h_{46}^{(4)} & h_{45}^{(4)} & h_{44}^{(4)} & h_{43}^{(4)} & h_{42}^{(4)} & K_{4,4}
\end{array}\right]=
$$

$$
=\left[\begin{array}{ccccccccc|c}
54.61 & 336.78 & 902.64 & 1378.23 & 1319.82 & 821.475 & 332.23 & 84.196 & 12.132 & 0.7585 \\ \hline
1.72 & 13.440 & 44.90 & 84.4200 & 98.5875 & 74.205 & 36.005 & 10.860 & 1.8475 & 0.1350 \\
0.300 & 2.350 & 7.875 & 14.8625 & 17.4375 & 13.200 & 6.4500 & 1.9625 & 0.3375 & 0.0250 \\ \hline
-100.33 & -758.22 & -2433.0 & -4362.15 & -4820.22 & -3405.6 & -1539.1 & -429.6 & -67.29 & -4.518 \\
-38.00 & -288.33 & -929.55 & -1675.22 & -1861.33 & -1322.3 & -600.67 & -168.3 & -26.44 & -1.777 \\
-12.00 & -92.00 & -300.33 & -549.333 & -621.00 & -450.00 & -209.00 & -60.00 & -9.666 & -0.666 \\ \hline
45.00 & 408.00 & 1485.50 & 2899.50 & 3401.81 & 2509.87 & 1170.86 & 334.50 & 53.31 & 3.6250 \\
-85.50 & -612.75 & -1872.8 & -3222.06 & -3439.1 & -2358.7 & -1039.0 & -283.6 & -43.56 & -2.875 \\
18.00 & 147.00 & 501.50 & 938.250 & 1064.00 & 761.50 & 345.50 & 96.250 & 15.00 & 1.000 \\
-18.00 & -129.00 & -390.50 & -658.750 & -684.75 & -456.00 & -195.00 & -51.75 & -7.750 & -0.500
\end{array}\right]
$$

To complete the example, let us calculate the determinant. We get it from the formula (6):

$$
\det(V)=\prod_{1\le i<j\le r}\left(\lambda_j-\lambda_i\right)^{n_i n_j}=\prod_{j=2}^{r=4}\prod_{i=1}^{j-1}\left(\lambda_j-\lambda_i\right)^{n_i n_j}
$$

Thus:

$$
\begin{aligned}
\det(V)&=\left[\prod_{i=1}^{1}\left(\lambda_2-\lambda_i\right)^{n_i n_2}\right]\cdot\left[\prod_{i=1}^{2}\left(\lambda_3-\lambda_i\right)^{n_i n_3}\right]\cdot\left[\prod_{i=1}^{3}\left(\lambda_4-\lambda_i\right)^{n_i n_4}\right]= \\
&=\left[\left(\lambda_2-\lambda_1\right)^{n_1 n_2}\right]\cdot \\
&\quad\cdot\left[\left(\lambda_3-\lambda_1\right)^{n_1 n_3}\left(\lambda_3-\lambda_2\right)^{n_2 n_3}\right]\cdot \\
&\quad\cdot\left[\left(\lambda_4-\lambda_1\right)^{n_1 n_4}\left(\lambda_4-\lambda_2\right)^{n_2 n_4}\left(\lambda_4-\lambda_3\right)^{n_3 n_4}\right]=-337.5
\end{aligned}
$$